\documentclass{scrartcl}
\usepackage[utf8]{inputenc}
\usepackage[T1]{fontenc}
\usepackage{graphicx}
\usepackage{longtable}
\usepackage{wrapfig}
\usepackage{rotating}
\usepackage[normalem]{ulem}
\usepackage{amsmath}
\usepackage{amssymb}
\usepackage{capt-of}
\usepackage{hyperref}
\usepackage{xcolor}
\usepackage{tikz}
\usepackage{soul}
\usepackage{listings}
\usepackage[version=3]{mhchem}
\usepackage{doi}
\usepackage{amsmath}
\usepackage[british, ]{babel}

\usepackage[british]{babel}
\usepackage[backgroundcolor=yellow!10!white]{mdframed}
\lstdefinelanguage{Julia}%
{morekeywords={@enum, abstract,Array,break,%
case,catch,const,continue,do,else,elseif,%
end,export,false,Real,for,function,%
immutable,import,importall,if,in,%
macro,module,mutable,new,%
otherwise,quote,Point,
print,return,%
struct,switch,true,try,typealias,%
using, Vector, while},%
aboveskip=12pt,%
basicstyle=\ttfamily\small,%
belowskip=12pt,%
breaklines=true,%
captionpos=t,%
commentstyle=\color{blue!50!black},%
emph={Agent, MessageType, Message@enum,%
},%
emphstyle=\color{red!50!black},%
firstnumber=1,
frame=lines,%
numbers=left,%
sensitive=true,%
alsoother={$},%
keywordstyle=\color{green!50!black},%
morecomment=[l]\#,%
morecomment=[n]{\#=}{=\#},%
morestring=[s]{"}{"},%
morestring=[m]{'}{'},%
stringstyle=\color{blue!50!black},%
}[keywords,comments,strings]
\author{Eric S Fraga\thanks{Corresponding author: \url{e.fraga@ucl.ac.uk}} \and Veerawat Udomvorakulchai \and Miguel Pineda \and Lazaros G. Papageorgiou \and \\ Sargent Centre for Process Systems Engineering \\ Department of Chemical Engineering \\ University College London (UCL)}
\date{}
\title{A multi-agent system for hybrid optimization}
\hypersetup{
 pdfauthor={Eric S Fraga \and Veerawat Udomvorakulchai \and Miguel Pineda \and Lazaros G. Papageorgiou \and \\ Sargent Centre for Process Systems Engineering \\ Department of Chemical Engineering \\ University College London (UCL)},
 pdftitle={A multi-agent system for hybrid optimization},
 pdfkeywords={},
 pdfsubject={},
 pdfcreator={Emacs 31.0.50 (Org mode 9.7.17)}, 
 pdflang={En_Gb}}
\begin{document}

\maketitle
\begin{abstract}
Optimization problems in process engineering, including design and operation, can often pose challenges to many solvers: multi-modal, non-smooth, and discontinuous models often with large computational requirements.  In such cases, the optimization problem is often treated as a \emph{black box} in which only the value of the objective function is required, sometimes with some indication of the measure of the violation of the constraints.  Such problems have traditionally been tackled through the use of \emph{direct search} and \emph{meta-heuristic} methods.  The challenge, then, is to determine which of these methods or combination of methods should be considered to make most effective use of finite computational resources.

This paper presents a multi-agent system for optimization which enables a set of solvers to be applied simultaneously to an optimization problem, including different instantiations of any solver.  The evaluation of the optimization problem model is controlled by a scheduler agent which facilitates cooperation and competition between optimization methods.  The architecture and implementation of the agent system is described in detail, including the solver, model evaluation, and scheduler agents.  A suite of direct search and meta-heuristic methods has been developed for use with this system.  Case studies from process systems engineering applications are presented and the results show the potential benefits of automated cooperation between different optimization solvers and motivates the implementation of competition between solvers. 

\end{abstract}
\section{Introduction}
\label{sec:org03961d0}

Optimization problems in process engineering, including design and operation, can often pose challenges to many solvers.  In particular, these problems may be multi-modal, non-smooth, and even discontinuous.  Identifying the best optimization method for any given problem is difficult if not impossible.  For a given problem, typically, methods will be tried until one is found to work.  Further, for some problems, it has been found that hybrid approaches are necessary to obtain the best solutions.

Hybrid methods bring together more than one method to work together to solve the problem.  One approach is to use two or more methods sequentially.  Examples of this include the use of a graphical interface with an embedded simulated annealing method to identify good initial guesses for a subsequent mathematical programming approaches, such as DICOPT or SBB with different NLP solvers such as MINOS or CONOPT \cite{fraga-papageorgiou-2006a}.  Another example is the design of heat exchanger networks using a variety of \emph{direct search} and \emph{metaheuristic} methods as well as user interaction to guide the search \cite{fraga-2006a,fraga-zilinskas-2002a}.

For problems which are non-smooth, including those with discontinuities, the evaluation of the objective function (or functions) and constraints is often as a \emph{black box}.  In this case, given a point in the search domain, \(d\), the objective function and constraints are evaluated the values determined.  No gradient information is explicitly available.  Such problems are traditionally solved using direct search and meta-heuristic methods where gradient information is not required.  Mathematical programming approaches can also be applied through the use of numerical differentiation to approximate the gradients.

This paper addresses the challenge of identifying the best optimization method or combination of methods for a given optimization problem.  A multi-agent system has been implemented which allows optimization solvers to work together to attempt to find the best solution to the given problem.  The design of the agent system aims to be generic but limited to methods which are suitable for black box problem specifications.

In what follows, the optimization problem definition is

\begin{align}
\label{optimizationproblem}
\min_{d \in \cal{D}} \, & z = f(d) \\
& g(d) \le 0 \nonumber
\end{align}

\noindent where we have assumed that \(g(d)\) includes also equality constraints for consistency with the implementation described below.  There are no assumptions on the properties of the search space, \(\cal{D}\): there are no assumptions related to the continuity, convexity, and smoothness of this space.  The decision variables, \(d\), may be real and/or integer values.  Finally, there are also no assumptions about the properties of the model for the problem, comprising both the objective function, \(f(d)\), and the constraints, \(g(d)\).  In the examples below, we consider problems with both real and integer variables and also with differential equations for the evaluation of the model.
\section{A multi-agent system}
\label{sec:org0b215c9}
A multi-agent system is software that enables autonomous pieces of software, known as \emph{agents}, to interact through the sending and receiving of \emph{messages} \cite{nwana1996software,bradshaw1997introduction}.  Multi-agent systems have been developed for engineering design problems, see for instance \cite{Zhang2020,10.1115/1.3013847,SIIROLA2003334,SIIROLA20031801} and references cited therein.  These previous implementations have agents working together to solve a design problem by either decomposing the search domain, by individual agents exploring the local domain around a point, or different agents addressing different regions of a Pareto frontier for multi-objective optimization problems.  Together, the agents act as a single optimization method.  This can be an effective approach, especially with multi-core systems where each agent is able to execute simultaneously with other agents.  However, such approaches presume that the overall search method is appropriate for the problem.  In other cases \cite{JULKA20021771,JULKA20021755,MELE2007722}, agents are used to model the problem, where each agent implements the behaviour of, for instance, different processing steps.

The aim of the multi-agent system proposed here is different: it is to enable different solvers and, in particular, different types of solvers to work independently but also cooperatively to solve a given design problem.  This is similar to an existing agent based system \cite{GEBRESLASSIE2017194} but allowing for heterogeneous optimization agents.  In the longer term, the design of the multi-agent system presented here is meant to support solvers competing for the computational resource, in the case where the evaluation of the objective function (and constraints) may be computationally expensive.  See \cite{pineda-etal-2021a,beck-etal-2015a} for examples of computationally expensive problems.  A secondary aim and bonus of the proposed multi-agent system is to minimise the changes required of any given optimization method to work within the agent system.  

The key design decision for the multi-agent system is to decouple the evaluation of both the objective function and the constraints from the solvers.  Normally, for a black box optimization method, the function implementing the objective and the constraints would be passed directly to the solver.  In the agent system, the solvers will be given a proxy function which interacts with the agent system to evaluate points in the search space.

\begin{figure}[hbtp]
\centering
\includegraphics[width=0.9\linewidth]{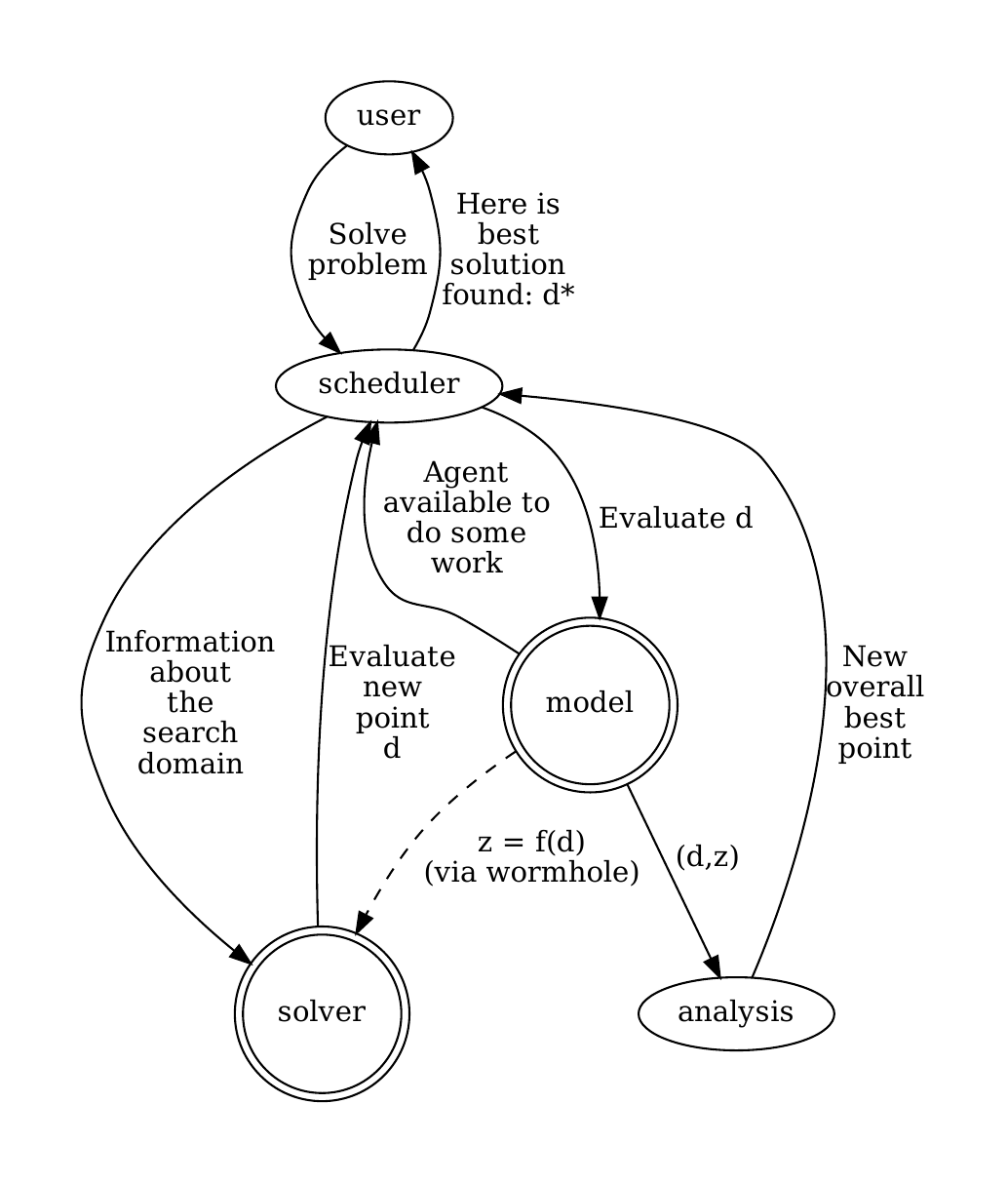}
\caption{\label{absdiagram}The set of agents in the agent based system including the scheduler, the model evaluation agents, the solver agents, and the analysis agent, showing the communications links between the agents.  Solid lines are persistent communication links; the dashed line represents ephemeral links to receive the result of the evaluation of the model, \emph{f(d)}, at a given point, \emph{d}, in the search domain.  The best solution found is indicated by \emph{d\textsuperscript{*}}.}
\end{figure}

The overall structure of the system is shown in Figure \ref{absdiagram}.  There are four types of agents along with the \emph{user} who defines the problem and invokes the system.  In the diagram, a single oval represents an agent (including the user) which is present as a single instance in the system.  A double oval represents an agent which can have one or more instances.  Solid lines represent persistent communication channels between the agents with the dashed line representing an ephemeral link between the model evaluation and individual solvers.  The four types of agents are described in more detail below.  The user is an actor in this system but does not, at this stage of development of the multi-agent system, get involved in the solution of the problem although such interaction has been shown to be useful for some problems \cite{fraga-2006a}.

\begin{lstlisting}[aboveskip=12pt,basicstyle=\ttfamily\small,belowskip=12pt,breaklines=true,captionpos=t,numbers=left,firstnumber=last,frame=lines,language=julia,label=messagestructure,caption={The structure defining messages along with the enumeration of the types of messages in Julia code.},captionpos=b]
@enum MessageType begin
    ANALYSESOLUTION      # for the analysis agent
    EVALUATEPOINT        # wish to evaluate the point
    OBJECTIVEVALUE       # the result of model on point
    REQUESTPOINT         # model evaluator available
    RETRIEVEBEST         # ask for best from analysis agent
    SHAREBEST            # feature: send best to all solvers
    STATISTICSBEST       # sending best known
end

struct Message
    type::MessageType    # type of message
    from::Agent          # agent sending the message
    content              # the actual message content
end
\end{lstlisting}

A multi-agent system is based on messages passing between agents.  The message structure defined for this system is shown in Listing \ref{messagestructure} as implemented in Julia \cite{bezanson-etal-2017a}.  Each message has a \texttt{type}, a reference to the agent which sent the message, and the content of the message.  Julia supports \emph{dynamic typing} and full introspection is possible.  This enables us to define the content without a type and it is up to the receiving agent to make use of this content based on the message type.  This is straightforward to implement as Julia supports \emph{multiple dispatch}.  The use and need for the different types of message will be highlighted in the description of the individual agents.
\subsection{The solver agent}
\label{sec:org0090783}

A solver agent represents an individual optimization method with a particular combination of parameter settings, should the method have such settings.  Each solver will be invoked with an initial population of points in the search space, usually randomly generated.  Every solver will be invoked with the same initial population.  It is assumed that each optimization method will be passed the function which implements the objective function and the constraints.  The signature of this function is

\begin{verbatim}
(z, g) = f(d, parameters)
\end{verbatim}

\noindent
where \texttt{d} is the point in the search space to evaluate and \texttt{parameters} represents problem specific values that may be necessary for the evaluation of the model at any point in the search space.  This function is expected to return the value or values of the objective function and an indication of feasibility.  A point is considered feasible if the value of \texttt{g} is \(\le 0\) and infeasible otherwise.  The actual value of \texttt{g} may be used by an optimization method to guide the search, for instance to attempt to find a feasible point from an infeasible point.  The Fresa plant propagation algorithm \cite{fresa-review-2021}, for instance, uses the measure of infeasibility to rank solutions in the assignment of a fitness value used for the selection of solutions to propagate.

In the multi-agent system, the solver is given a function which is not the actual model for the problem but a proxy function.  This proxy function sends a message to the \emph{scheduler} agent (see below) with a type of \texttt{EVALUATEPOINT} and the content consisting of two pieces of information: the point in the search space to evaluate and an ephemeral \emph{channel} to communicate the result of the evaluation, \texttt{(z, g)}, back to this function which then returns this value back to the optimization method that invoked the function.  In Figure \ref{absdiagram}, this ephemeral channel is indicated with a dashed line and labelled as a \emph{wormhole}.

There will always be more than one solver agent as the aim is to explore the effect of interactions between the different solvers.  
\subsection{The model evaluation agent}
\label{sec:org8d766a1}

The model evaluation agent is invoked with the function that implements the objective function and the constraints (Eq. \ref{optimizationproblem}), along with any parameters that function may require.  The agent is also provided with the channel through which it may communicate with the analysis agent.

There may be several instances of this agent, each evaluating the same model.  This allows the agent system to make use of modern multi-core computers for concurrent evaluation of different points.

Each instance of the model evaluation agent performs a loop:
\begin{enumerate}
\item Send a message with type \texttt{REQUESTPOINT} to the scheduler, indicating that this model agent is ready to do more computational work.
\item Wait for a message from the scheduler.  This message should be of type \texttt{EVALUATEPOINT} which will have, as contents, two pieces of information, the point in the search space to evaluate and an ephemeral Julia channel to use to send the outcome of the model evaluation.
\item Evaluate the given point in the search space.
\item Send the outcome of the evaluation as a message of type \texttt{OBJECTIVEVALUE} to the solver agent, using the ephemeral channel provided in the request for evaluation message, and also to the analysis agent.
\item Go to step 1
\end{enumerate}
\subsection{The analysis agent}
\label{sec:org73b956f}

Each time a point in the search space is evaluated, the point and the result of the evaluation are sent to the analysis agent.  In the current version of the multi-agent system, the analysis agent simply checks to see if the point is \emph{better} than any found to date.  For single criterion problems, this is a simple check: the agent simply keeps track of the overall best solution found by any solver and updates it when a better point is encountered.  For multicriteria problems, the agent keeps a set of the non-dominated points found and the new point is considered \emph{better} if it is not dominated (and the set of non-dominated points is updated).  When a better point is found, it is sent to the scheduler.
\subsection{The scheduler agent}
\label{sec:orgfb76d72}

\begin{figure}[hbtp]
\centering
\includegraphics[width=0.8\textwidth]{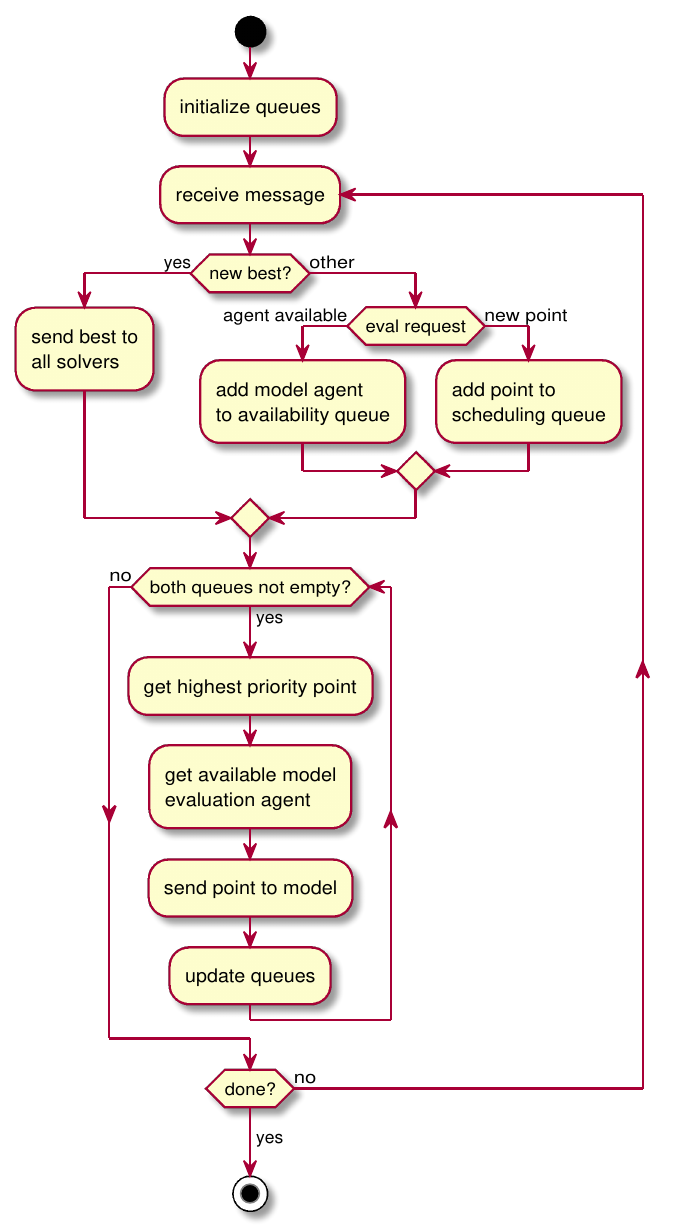}
\caption{\label{schedulerflowchart}The flow chart for the scheduler agent which assigns work to the model evaluation agents based on requests from the solvers and communicates with all solvers when new \emph{best} solutions are identified.}
\end{figure}

At the centre of the multi-agent system is the scheduler; see Figure \ref{schedulerflowchart} for a flowchart of this agent's behaviour.

The main purpose of this agent is to process requests from the solvers for the evaluation of points in the search domain.  When such a request is received, it is added to an evaluation queue.  This queue is used to allocate the evaluation of these points when model evaluation agents become available.  The availability of model evaluation agents is handled by a different queue.

The queue of model evaluation requests is priority based and consists of a vector of \emph{first-in first-out} (FIFO) queues.  Each solver is given a priority, a number between 1 and a maximum (currently 10).  When a request comes in, the request is added to the queue at that priority level.  When the schedule wishes to send a point to a model evaluation agent, the highest priority point is chosen.  The first point in each FIFO queue is then promoted to the next higher priority level if not already at the top level; this is the \emph{update queues} action in Figure \ref{schedulerflowchart}.  This approach to scheduling the requests is similar to that used to implement multi-tasking in operating systems \cite{madnick-donovan-1974}, ensuring that all solvers eventually have their model evaluation requests satisfied.  In the case studies presented below, all solvers have the same priority so the behaviour overall is as a FIFO queue but the implementation is more powerful than that.

The queue for the availability of the model agents is itself a simple FIFO queue.  The order in which the various instances of model evaluation agents are given points to evaluate has no effect on the overall evaluation of points, given that all agents are executing on the same computer.

The scheduler also enables sharing of information between solvers: it sends any new best solution communicated to it by the analysis agent to all the solvers.  It is up to the solvers individually to decide what to do with this information.  The individual behaviour of the solvers and how they use this information is described in the next section.
\subsection{Implementation}
\label{sec:org63d8152}

The system has been implemented in the Julia language \cite{bezanson-etal-2017a}, a high performance general purpose language with high level data structures and funtional programming features.  It supports multi-processing and multi-threading computation natively and provides all the elements required for the implementation of an agent based system.  Specifically, two key elements of the Julia system are used in the implementation:

\begin{description}
\item[{Threads}] parallel computation using light-weight tasks which can communicate with other threads and share memory between them, and
\item[{Channels}] which are simple first-in first-out queues of messages with any type of content that can be accessed by different threads.
\end{description}

\noindent
So, for instance, a channel can be created to allow the scheduler to send messages to a specific model agent and another channel created to send messages to the analysis agent.  As an example,

\begin{verbatim}
analysischannel = Channel{Message}(2*(nm+1))
\end{verbatim}

\noindent
creates a channel which can store a specified number of messages (two for each model agent plus one).  Any agent that has access to this channel, that is any that has this variable in \emph{scope}, can place messages into the channel, using the \texttt{put!} function, and can retrieve messages from the channel, using the \texttt{take!} function.

Then, a thread can be started by \emph{spawning} where a particular function, with its arguments, is specified.  The thread will start and the function will be executed.  The thread will terminate when the function terminates, which is why all agents, except the scheduler, loop forever.  Only the scheduler eventually terminates which then terminates the solution process.

So, for example, a model agent can be instantiated by

\begin{verbatim}
Threads.@spawn modelagent(model, schedulerchannel, analysischannel)
\end{verbatim}

\noindent
where \texttt{modelagent} is simply our own Julia function with three arguments: the \texttt{model}, a function which evaluates the actual objective function and associated constraints (see equation \ref{optimizationproblem}), \texttt{schedulerchannel} is the channel through which the scheduler will send requests for model evaluations, and \texttt{analysischannel} is the channel the model evaluation agent will use to send results to the analysis agent.

Multi-thread applications require careful implementation to avoid ending up in a deadlock situation where two or more processes, or threads, are waiting on messages from processes that are themselves waiting for messages \cite{Dijkstra1971}.  The algorithms described above for each of the agents have been designed to avoid such a situation arising.  Further, compared with other multi-agent systems \cite{GEBRESLASSIE2017194}, the system presented here has no global data and therefore avoids any need for synchronisation on the access to these data.  Each solver works autonomously and all data transfer is via the coordinated messages.
\section{The solvers}
\label{sec:org03e89c3}

Two types of optimization methods, or solvers, are considered: meta-heuristic methods inspired by nature \cite{fraga-book-2022} and direct search methods \cite{kelley-1999}.  Meta-heuristic methods are often good at exploring the search space globally but less effective at tuning good solutions.  Direct search methods often have the opposite charateristics: they may not be effective at global exploration but they may excel at fine tuning.  The hope is that by providing a range of solvers with different characteristics, challenging optimization problems can be solved without \emph{a priori} decisions about which solver may be best or without requiring manual intervention to, for instance, use one solver and then use the result of that solver to initialise a second solver.

The suite of methods has been implemented in the Julia language and the methods have been instrumented to enable the agent based system to share information when appropriate with each solver.  None of the implementations could be considered \emph{best in class} but they have all been implemented to provide a demonstration of the class of methods they represent.  They are sufficiently different to demonstrate the potential of the multi-agent system.

There may be more than one instance of each solver in the system.  Different instances of the same solver can have different values for key parameters that affect their performance.  For instance, the population size for any population based approach is often a parameter that needs to be chosen by the user and for which there is often little or no guidance.  Having different instances of a solver, each with a different population size, may avoid the need to guess the correct value of this parameter and may also lead to the avoidance of premature convergence during the search.

The full list of solvers used in this paper follows:

\begin{description}
\item[{genetic algorithm}] (GA) inspired by Darwinian evolution \cite{holland-1975} using \emph{crossover} and \emph{mutation} with fitness based \emph{selection} to explore the solution space and exploit the better solutions.  In what follows, the size of the population modelled is the key parameter that is explored by having multiple instances of this solver.  This method incorporates new \emph{best} solutions found by other solvers by simply adding the new solution to the population before the start of a generation.  This is similar to the metaphor of multi-island genetic algorithms \cite{1041554}.

\item[{plant propagation algorithm}] (PPA) inspired by the propagation of strawberry plants using \cite{salhi-fraga-2011a} \emph{runners} where the number of runners propagated by a solution is proportional to the fitness of that solution and the length of the runners inversely proportional to the fitness.  This combination leads to a balance between exploitation and exploration of the solution space.  The number of solutions to consider for propagation at each iteration is the key parameter explored.  As with the GA, new solutions sent by the scheduler are added to the population at the start of an iteration.  In terms of nature inspiration, this could be viewed as incorporating a further means of propagation used by plants: \emph{seeds} carried over from other areas by animals or the wind \cite{sulaiman-salhi-2016a}.

\item[{particle swarm optimization}] (PSO) based on the behaviour of, for instance, flocks of birds in the sky \cite{488968} where the evolution of each member of the population is based on each member having a velocity (both magnitude and direction) that is updated based on where each solution is relative to the best solution in the swarm and its own history in the evolution of the population.  Solutions from the schedule replace the currently worst solution in the swarm at the start of an iteration.

\item[{steepest descent}] (SD) is a gradient based hill-climbing method which has been implemented using finite differences for the numerical approximation of the gradient \cite{Burden1989}.  The initial population is treated as a \emph{multi-start} optimization problem where a search is started with each member of this population in turn.  New solutions from the scheduler are added to the list of starting points yet to be attempted and the list is treated as a last-in first-out queue so that focus switches to the currently known best area in the search space.  For single objective optimization problems, there will be only one instance of this solver.  However, for multi-objective problems, there will be several instances, each attempting the multi-objective problem as a single objective problem but with different weightings between the objectives.

\item[{coordinate search}] (CS) which searches along each dimension in turn, using a simple line search method (also used by the steepest descent method) \cite{kelley-1999}.  The initial population and any points from the scheduler are treated as for the steepest descent method.  As with the steeptest descent method, the number of instances of this solver will be 1 for single objective problems and more than 1 for multi-objective problems.
\end{description}
\section{Case studies}
\label{sec:org8c0b2d1}

Recall that one aim of the agent system is to enable the application of multiple methods for the solution of an optimization problem when there is no \emph{a priori} knowledge of which method may be most effective.  A second aim is to make most effective use of multiple methods through the sharing of information about the search space, dynamically while the search is in progress.  Two case studies are presented, one on the design of heat exchanger networks, a problem that combines both nonlinear and combinatorial aspects, and a second which is a computational expensive multi-objective design problem, relying on dynamic simulation as an integral part of the evaluation of the objectives.
\subsection{Heat exchanger network design}
\label{sec:org515c242}
Two case studies in the design of heat exchanger network synthesis are presented to explore the impact of the sharing of information to allow solvers to cooperate.  One is a small, 2 hot stream and 2 cold stream, example and the other is a larger problem with 5 cold and 5 hot streams based on case study 5 in \cite{fraga-2009a}.  The first is solved as an MINLP and the second as a pure NLP using a fully continuous representation of the superstructure \cite{fraga-2006a}.  Heat exchanger network design problems are combinatorial in nature and solutions often have a subset of design variables at the constraints.  Previous work \cite{fraga-2006a} has shown that manually applying a sequence of methods can lead to better final outcomes.  By using the multi-agent system, it will be interesting to see if and how the methods can cooperate to obtain good solutions.

The solvers considered are a subset of the methods described above, with multiple instances of the meta-heuristic methods with different parameter settings.  The base population size, \texttt{np}, is the number of design variables for both problems.  This is adjusted for different instances of the solvers: for the GA are \texttt{np} (small) and \texttt{5 × np} (large); for the PPA, the sizes are \texttt{np/2} (small) and \texttt{2 × np} (large), where these numbers actually define the number of solutions to propagate at each iteration, not the population size.  There is only one instance of each of the direct search methods.  All methods are started with the same randomly generated population of size \texttt{np}.

The stopping criterion for the search is 60 thousand messages received by the scheduler agent, which can otherwise be thought of as 60 thousand iterations over the main loop in the scheduler flowchart (Figure \ref{schedulerflowchart}).  This equates to on the order of 5 thousand function evaluations by each of the 6 solvers.

\begin{figure}[hbtp]
\centering
\includegraphics[width=0.9\linewidth]{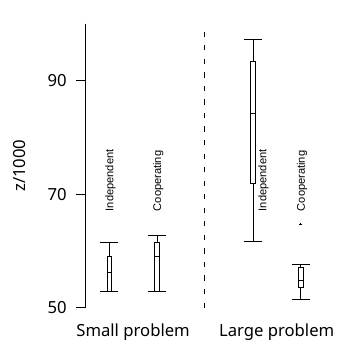}
\caption{\label{boxplotfigure}Box plot of the outcomes of 10 runs for both problems, comparing the outcomes for both when solvers are fully independent and when they are cooperating by sharing improved solutions when found.}
\end{figure}

Figure \ref{boxplotfigure} presents the box plots of the variation of the objective function value for the best solution found over 10 runs, both when each solver works independently or when the solvers share any new better solution found.  For the smaller problem, not sharing information actually leads to better outcomes, although only marginally.  It would seem that sharing may lead to premature convergence in some cases.  However, for the larger problem, sharing information leads to better solutions more quickly.

\begin{figure}[hbtp]
\centering
\includegraphics[width=0.9\linewidth]{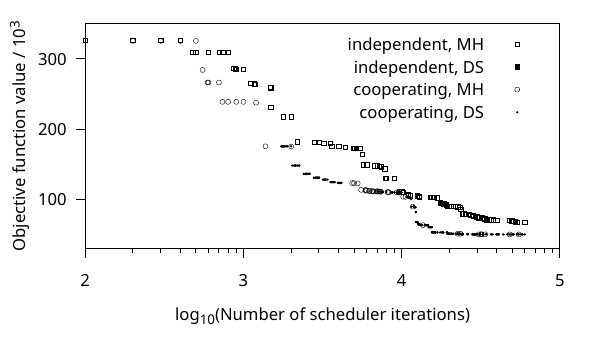}
\caption{\label{10spevolutionfigure}Evolution of best solution found as a function of the number of scheduler iterations indicating the solver responsible for the identification of the best solution to that point. The top trail of points is for when the solvers operate independently; the bottom trail is the outcome when the best solutions found are shared.  MH is meta-heuristic and DS is direct search.}
\end{figure}

Figure \ref{10spevolutionfigure} shows the evolution of the best solution found as a function of the number of scheduler iterations for both independent and cooperating cases, indicating which particular type of method (meta-heuristic or direct search) is responsible for each improvement in the objective function value of the best solution found at that point.

For the independent case, only the meta-heuristic methods are able to find good solutions.  When the methods cooperate, both meta-heuristic and direct search methods are able to find improved solutions.  It appears that when one of the meta-heuristic methods identifies a better solution, this solution is immediately improved by one of the direct search solvers.  The meta-heuristic methods are good at \emph{exploration} and direct search methods are good at \emph{exploitation}.  Further, we see that the cooperating approach leads to not only a better final solution, these solutions are obtained in approximately 20-30\% of the overall computational time compared with not cooperating.  The automatic hybridization of the methods through cooperation is both effective and efficient.

It should be noted that, as some of the methods are stochastic and the initial population is also randomly generated, repeated attempts will lead to different outcomes.  However, the outcomes presented here are representative.

For the larger case study, the behaviour when solutions are shared is what we expected.  The direct search methods have better exploitation properties.  It should be noted, however, that the results obtained in this case study are not as good as the known optimum for this problem.  This is most likely due to the rather generic untuned nature of the different implementations.  It is expected that if the individual solvers in the multi-agent system were better implementations (e.g. NSGA-II \cite{deb-2000}, Fresa \cite{fresa-review-2021}, BFGS \cite{fletcher-1987a}, other direct search methods \cite{kelley-1999}), the outcomes would be better.
\subsection{Design of a micro-analytic system}
\label{sec:org05c4a75}

A model-based optimization approach for the design of microfluidic separation processes together with an online detection tool such as thermal lens microscopy was previously developed \cite{pineda-etal-2021a}.  The application of the resulting micro-analytic analysis system was for the processing of waste from the decommissioning the Fukushima Dai-chi Nuclear Power Plants where there is a need to detect the presence of radionuclides in water.

Microfluidic systems are based on channels with characteristic dimensions measured in micrometres.  In these channels, transport resistances are reduced and heat and mass transfer rates intensified.  Flow patterns can be established that increase the transport rates and interfacial area while facilitating the separation of phases at the end of the microchannels. The well defined flow patterns and the laminar flow conditions allow precise modelling and characterisation of such units.  The models are based on differential equations, requiring significant computational effort for the evaluation of individual designs.

Overall, the model is both differential and stochastic, with the optimization problem considered as an example of \emph{design under uncertainty} \cite{sahinidis-2004a,grossmann-etal-2016a}.  There are two objectives for the design: minimise the number of false readings and maximise the detection of radionuclides (which is equivalent to finding the minimum concentration detected).  False readings may be positive readings, which means that the detector identifies the presence of radionuclides when none is present, or negative readings, where radionuclides present but not detected by the device.  The latter are more serious than the former.  Positive readings may lead to more expense (e.g. more detailed analysis) whereas negative readings may give a dangerous sense of false security.

The optimization problem is an example of where the multi-agent system may be most appropriate: unlike many other problems, including the heat exchanger network design problem above, where many function evaluations are possible, solving this problem will need to be achieved using a relatively small number of function evaluations.  In this case, it is desirable for each function evaluation to be as effective as possible.  A further challenge is the need to handle multiple objectives, leading to the need for larger populations for the meta-heuristic methods and multiple instances of the direct search methods to tackle a single objective defined by a weighted sum of the individual objectives, using different weights in each instance.  Table \ref{mutassizeparametertable} summarises the population size for meta-heuristic methods, explicitly as a function of n\textsubscript{p} (chosen to be 20 for this problem) for the genetic algorithm and the particle swarm optimization methods and implicitly via the setting of the number of solutions to propagate for the plant propagation algorithm. 

\begin{table}[hbtp]
\caption{\label{mutassizeparametertable}The population size, for genetic algorithm (GA) and particle swarm optimization (PSO) methods, and number of solutions to propagate for the plant propagation method for a given value of n\textsubscript{p}.}
\centering
\begin{tabular}{lr}
Solver & Size parameter\\
\hline
small GA & 2 n\textsubscript{p}\\
big GA & 5 n\textsubscript{p}\\
small PPA & 0.5 n\textsubscript{p}\\
big PPA & 2 n\textsubscript{p}\\
small PSO & n\textsubscript{p}\\
big PSO & 5 n\textsubscript{p}\\
\end{tabular}
\end{table}
\subsubsection{Measures of quality of non-dominated set}
\label{sec:orgfc59306}
It is straightforward to compare the outcome of different methods or approaches for single objective optimization problems: simply compare on the basis of the single scalar value of the objective function for the best solution identified in each attempt.  It is less straightforward to compare outcomes for a multi-objective problem.  The \emph{best} solution for a multi-objective problem is not a single scalar value but a vector of values, known as the \emph{set of non-dominated points}.

Comparing such vectors of outcomes can be done using a number of measures from the literature:

\begin{description}
\item[{Hypervolume}] (HV) \cite{1583625,ZitzlerandThieleMOEA}  is a metric used to quantify the volume of the objective space that is dominated by a given set of solutions. It encapsulates both the proximity of solutions to the true Pareto front and, to some extent, the diversity of solutions across the objective space, within a single scalar value \cite{7360024}. However, when the true Pareto front is unknown, its application requires preference information to select an appropriate reference point, which may result in misleading values. It is worth noting that whenever one solution set completely dominates another, the hypervolume of the former will be higher than that of the latter; however, this only holds when the true Pareto front is known.

\item[{Area}] This metric calculates the area of a 2D polygonal shape projected onto the axis of one of the objectives by applying the trapezoidal rule to the points in a given set of solutions. It quantifies how much of the bi-objective space is covered by the non-dominated solutions. A larger area implies that the solutions are spread and explore a broader range of trade-offs. However, this metric does not necessarily measure how close the points are to the true Pareto front. It is also sensitive to outlier solutions, which may result in misleading values of the computed area.

\item[{Average Distance}] This metric computes the average of the Euclidean distances to the reference (utopia) point. This measure indicates how close the solutions are to the reference point. However, this metric does not provide any insights into the diversity of the solution within the objective space and is sensitive to the location of the utopia point.

\item[{Generational Distance}] (GD) \cite{van1998evolutionary} is a metric used to measure how far a given set of solutions is from the true Pareto front. This metric provides a single scalar value by computing the average Euclidean distance from each point in the solution set to its nearest point on the true Pareto front. A smaller value of GD is desirable, as it indicates better convergence toward the Pareto front. However, a key limitation of this metric is that it requires the knowledge of the true Pareto front of the problem and does not offer any information about the diversity of the solutions.
\end{description}
\subsubsection{Results}
\label{sec:orgcaa92f2}
Given the computational expense, where the evaluation of a single design can take anything from a second to hundreds of seconds depending on the representation of the stochasticity of the feeds considered, only 1000 function evaluations in total are allowed.  There are two instances of each meta-heuristic, with different population sizes; there are four instances of each direct search method, with weightings \(\omega\)=\{0.2, 0.4, 0.6, 0.8\} where the objective function used for the direct search methods is
\[ z = \omega z_1 + (1-\omega) z_2 \]
and where \(z_1\) is the sensitivity of the device, \(\phi\), a factor relative to the sensitivity of the actual sensing technology, and
\[z_2 = n_f + \frac{n_p}{1000}\]
with \(n_f\) being the number of false negatives and \(n_p\) the number of false positives.  The use of multiple instances of the direct search methods, with different weightings, is similar to previous agent based approaches for the identification of a discrete set approximating the Pareto frontier \cite{SIIROLA2003334}.

The choice of the second objective function highlights the difference in impact in having negative outcomes.  False positive readings occur when the analytic device would indicate that radionuclides were present in sufficient quantities to be dangerous when not actually present in those quantities.  This may lead to greater expense in dealing with the potentially radioactive material.  False negative outcomes are when the device would indicate that the amount of radionuclides present is below the threshold where in fact there is a greater amount actually present.  Such a negative reading may lead to dangerous situations where material is considered safe when it is not so.

\begin{table}[hbtp]
\caption{\label{mutasparetomeasurestable}Different measures of the set of non-dominated points in the final outcome of the multi-agent system for the design of the analytic system.  In all cases, except for the number of points, a lower value is better.  A larger value for the number of points may indicate a better representation of the set of non-dominated points with respect to the actual Pareto set for the optimization problem.}
\centering
\begin{tabular}{lrr}
Measure & Independent & Cooperating\\
\hline
hypervolume & 70.9 & 37.7\\
area & 1.5 & 1.1\\
average distance & 11.0 & 10.3\\
generational distance & 8.2 & 6.5\\
\hline
non-dominated points & 9.0 & 9.3\\
\end{tabular}
\end{table}

Table \ref{mutasparetomeasurestable} shows the difference in performance of the multi-agent system without sharing (\emph{independent}) and with sharing (\emph{cooperating}) of solutions during the evolution of the search.  For all measures, the outcomes with sharing of solutions are better than without.  However, these results are only an indication as the measures have generally been defined for outcomes where the global optimum Pareto frontier or set is known.  We have used the utopia point for these measures.

\begin{figure}[hbtp]
\centering
\includegraphics[width=0.9\linewidth]{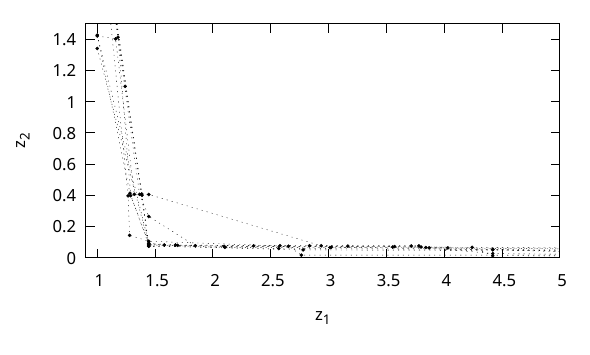}
\caption{\label{mutasnotsharingstochasticoutcomesfigure}Final sets of non-dominated points for 10 runs of the multi-agent system without sharing solutions between solvers.}
\end{figure}

\begin{figure}[hbtp]
\centering
\includegraphics[width=0.9\linewidth]{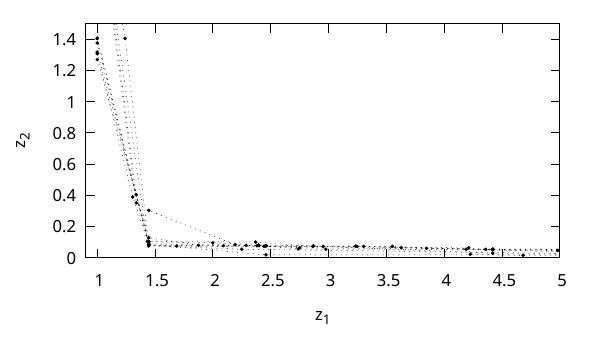}
\caption{\label{mutassharingstochasticoutcomesfigure}Final sets of non-dominated points for 10 runs of the multi-agent system sharing solutions between solvers.}
\end{figure}

A comparison can also be done graphically.  Figures \ref{mutasnotsharingstochasticoutcomesfigure} and \ref{mutassharingstochasticoutcomesfigure} present the final set of non-dominated points in 10 runs, first without sharing solutions between solvers and secondly with sharing of solutions.  The stochastic nature of the methods is evident in the distribution of possible outcomes.  However, we see that the use of sharing leads to more consistent outcomes and outcomes consistent with the measures shown in Table \ref{mutasparetomeasurestable}.

When looking in detail at the evolution of the search, we note the following:
\begin{itemize}
\item All methods other than the coordinate search method contribute to the search by finding new \emph{better} solutions.
\item When not sharing solutions, the vast majority of the solutions found are identified by meta-heuristic methods.
\item When solutions are shared, the direct search methods, specifically the steepest descent method, is able to find good points and contribute effectively to the search.
\item The steepest descent method contributes mostly towards the end of the search, improving on solutions obtained by the meta-heuristic methods.
\end{itemize}
\section{Conclusions}
\label{sec:org074e7dd}
Optimization problems that arise in the process industries may often exhibit properties that preclude the use of mathematical programming, properties such as nonconvexity, nonlinearity, discontinuities, and differential equations.  For those problems, a variety of meta-heuristic and direct search methods, that may treat the objective function as a black box, can be considered.  However, there is a large choice of these methods, each of which may also have a large number of tunable parameters.  As a result, choosing the appropriate method may be difficult \emph{a priori}.  Furthermore, each of these methods may have beneficial behaviour in different parts of the search space.

These two factors motivate the development of a multi-agent system to enable a variety of methods, with different instances of many, to be applied to a problem simultaneously. The simultaneous application of multiple methods and instances of methods leads naturally to the idea of sharing information while the search progresses.  This paper has shown that doing so leads to improved outcomes.

In the next steps, we wish to consider more complex scheduling algorithms in the scheduler agent and also explore what may happen when an element of competition is introduced where the solvers have to compete for the computational resource.  This is motivated by noting that the direct search methods help once good solutions are identified by the meta-heuristic methods so allocating equal computational resource to all methods all the time may not be the most effective approach.  A more adaptive approach may lead to more efficient use of the computational resource.
\subsection{Acknowledgements}
\label{sec:org6dc7481}
Mr. Veerawat is grateful for the support provided by the Royal Thai Government Scholarship, which covers his tuition fees and living stipends throughout the duration of his studies at University College London (UCL).  The authors also gratefully acknowledge funding received from UK Research \& Innovation through the grant EP/R019223/1 which supported the development of the models for the second case study.  

\bibliographystyle{plain}
\bibliography{/home/ucecesf/s/share/texmf/bibtex/bib/bibliography}
\end{document}